\documentclass{birkjour_t2}
\usepackage{float}
\usepackage{amsmath}
\usepackage{amsfonts}
\usepackage{amssymb}
\usepackage[english]{babel}
\usepackage{color}
\usepackage{graphicx}
\usepackage{tkz-euclide}
\usetkzobj{all}
\tikzset{elegant/.style={smooth,thick,samples=50,cyan}}
\tikzset{eaxis/.style={->,>=stealth}}
\usepackage{color,tikz}
\usetikzlibrary{decorations.pathreplacing,calc}
\usepackage{rotating}
\tikzset{liltext/.style={font=\tiny}}
\usepackage{lmodern}
\usepackage{amsthm}
\usepackage{graphicx}
\usepackage{mathrsfs}
\usepackage{makecell}
\usepackage{microtype}
\usepackage{enumerate}
\usepackage{mathscinet}
\usepackage{array}
\usepackage{multirow}
\usepackage{enumerate}
\usepackage[cal=boondoxo,bb=ams]{mathalfa}
\usepackage{hyperref}
\hypersetup{hidelinks}

\renewcommand{\theequation}{\arabic{section}.\arabic{equation}}

\newtheorem{thm}{Theorem}[section]

\newtheorem{lem}{Lemma}[section]

\newtheorem{rem}{Remark}[section]

\newcommand{\ml}{\mathcal}
\newcommand{\mb}{\mathbb}

\begin{document}

%
%
%
%
%
%
%
%

\title[Some remarks on blow-up of solutions to Nakao's problem]
 {Some remarks on blow-up of solutions to Nakao's problem}

%
%
\author[W. Chen]{Wenhui Chen}
\address{Institute of Applied Analysis, Faculty for Mathematics and Computer Science\\
	 Technical University Bergakademie Freiberg\\
	  Pr\"{u}ferstra{\ss}e 9\\
	   09596 Freiberg\\
	    Germany}
\email{wenhui.chen@student.tu-freiberg.de}
\author[M. Reissig]{Michael Reissig}
\address{Institute of Applied Analysis, Faculty for Mathematics and Computer Science\\
	Technical University Bergakademie Freiberg\\
	Pr\"{u}ferstra{\ss}e 9\\
	09596 Freiberg\\
	Germany}
\email{reissig@math.tu-freiberg.de}


%
\begin{abstract}
In this note we try to understand the blow-up of solutions to Nakao's problem by using nonlinear ordinary differential inequalities.
\end{abstract}

\maketitle
\section{Introduction}
The goal of this note is to understand nonexistence results for global (in time) solutions to the Cauchy problem for the weakly coupled system of semilinear damped wave equation and semilinear wave equation, that is, blow-up of solutions to the following weakly coupled system:
\begin{equation}\label{Nakao's problem}
\left\{
\begin{aligned}
&u_{tt}-\Delta u+u_t=|v|^p,&&x\in\mb{R}^n,\,\,t>0,\\
&v_{tt}-\Delta v=|u|^q,&&x\in\mb{R}^n,\,\,t>0,\\
&(u,u_t,v,v_t)(0,x)=(u_0,u_1,v_0,v_1)(x),&&x\in\mb{R}^n,
\end{aligned}
\right.
\end{equation}
where $n\geq1$ and $p,q>1$.

The problem of critical curve for the weakly coupled system \eqref{Nakao's problem} was proposed by Professor Misuhiro Nakao, Emeritus of Kyushu University (see \cite{NishiharaWakasugi2015,Wakasugi2017}). Here ``critical curve" means that the threshold condition of a pair of exponents $(p,q)$ for global (in time) existence of small data Sobolev solutions and blow-up of classes of local (in time) Sobolev solutions. Recently, by applying the test function method (e. g. \cite{MP01book,Zhang01}) the author of \cite{Wakasugi2017} proved that if the condition
\begin{align}\label{Condition Wakasugi}
\alpha_{\text{N,W}}:=\max\left\{\frac{q/2+1}{pq-1}+\frac{1}{2};\frac{q+1}{pq-1};\frac{p+1}{pq-1}\right\}\geq \frac{n}{2},
\end{align}
holds for spatial dimensions $n\geq1$, then local (in time) Sobolev solutions $(u,v)$ to \eqref{Nakao's problem} in general blow up in finite time. However, in the $p-q$ plane for a pair of exponents $(p,q)$ the curve $\alpha_{\text{N,W}}=\frac{n}{2}$ is optimal only when $n=1$
because the condition $\alpha_{\text{N,W}}\geq\frac{1}{2}$ means that every local (in time) Sobolev solution blows up for $1<p,q<\infty$. When $n\geq 2$, the condition \eqref{Condition Wakasugi} seems not to be optimal. Thus, if a pair of exponents $(p,q)$ does not satisfy the condition \eqref{Condition Wakasugi}, the questions of global (in time) existence or nonexistence of Sobolev solutions to Nakao's problem \eqref{Nakao's problem} are still open.

We sketch now some historical background to \eqref{Nakao's problem}. Since Nakao's problem \eqref{Nakao's problem} is related to a weakly coupled system  of semilinear damped wave equation and semilinear wave equation, we recall some results for weakly coupled systems for semilinear wave equations and semilinear damped wave equations, respectively, in the following.

On one hand, the following weakly coupled system of semilinear wave equations
\begin{equation}\label{weaklycoupledwave}
\left\{\begin{aligned}
&u_{tt}-\Delta u=|v|^{p},&\quad &x\in\mathbb{R}^n,\,\,t>0,\\
&v_{tt}-\Delta v=|u|^{q},&\quad &x\in\mathbb{R}^n,\,\,t>0,\\
&(u,u_t,v,v_t)(0,x)=(u_0,u_1,v_0,v_1)(x),&\quad &x\in\mathbb{R}^n,
\end{aligned}\right.
\end{equation}
for $n\geq1$ with $p,q>1$, has been widely studied in recent years. The papers \cite{DelS97,DGM,DM,AKT00,KT03,Kur05,GTZ06,KTW12} investigated that the critical curve for \eqref{weaklycoupledwave} is described by the condition
\begin{equation}\label{Critical curve wave}
\alpha_{\text{W}}:=\max\left\{\frac{p+2+q^{-1}}{pq-1};\frac{q+2+p^{-1}}{pq-1}\right\}=\frac{n-1}{2}.
\end{equation}
In other words, if $\alpha_{\text{W}}<\frac{n-1}{2}$, then there exists a unique global (in time) Sobolev solution for small data; else if $\alpha_{\text{w}}\geq\frac{n-1}{2}$, in general, local (in time) Sobolev solutions blow up in finite time.

On the other hand, let us consider the weakly coupled system of semilinear classical damped wave equations
\begin{equation}\label{weaklycoupleddamped}
\left\{\begin{aligned}
&u_{tt}-\Delta u+u_t=|v|^{p},&\quad &x\in\mathbb{R}^n,\,\,t>0,\\
&v_{tt}-\Delta v+v_t=|u|^{q},&\quad &x\in\mathbb{R}^n,\,\,t>0,\\
&(u,u_t,v,v_t)(0,x)=(u_0,u_1,v_0,v_1)(x),&\quad &x\in\mathbb{R}^n,
\end{aligned}\right.
\end{equation}
for $n\geq1$ with $p,q>1$. The critical curve for \eqref{weaklycoupleddamped} is described by the condition
\begin{equation}\label{Critical curve damped wave}
\alpha_{\text{DW}}:=
\max\left\{\frac{p+1}{pq-1},\frac{q+1}{pq-1}\right\}=\frac{n}{2}.
\end{equation}
which has been investigated by the authors of \cite{SunWang2007,Narazaki2009,Nishi12,NishiharaWakasugi}.

From the above results of critical curves for \eqref{weaklycoupledwave} and \eqref{weaklycoupleddamped}, we may expect that the critical curve for \eqref{Nakao's problem} is between \eqref{Critical curve wave} and \eqref{Critical curve damped wave}. However, we should underline that the critical curve for \eqref{Nakao's problem} is not a simple combination of \eqref{Critical curve wave} and \eqref{Critical curve damped wave} because the critical curve to \eqref{Nakao's problem} seems to be influenced by varying degrees between semilinear wave equation and semilinear damped wave equation.

Let us explain the difficulties to derive blow-up of solutions. To obtain blow-up of solutions to \eqref{weaklycoupledwave} when $\alpha_{\text{W}}\geq\frac{n-1}{2}$, the authors of \cite{DGM,DM,DelS97,KTW12} mainly applied generalized Kato's lemmas to a system of nonlinear ordinary differential inequalities and constructed some contradictions. However, the authors of \cite{NishiharaWakasugi,SunWang2007} derive blow-up of solutions to \eqref{weaklycoupleddamped} when $\alpha_{\text{DW}}\geq\frac{n}{2}$ by applying the test function method. So, it is interesting to see what method is suitable for us to derive blow-up of solution to \eqref{Nakao's problem}. Although \cite{Wakasugi2017} has applied the test function method to \eqref{Nakao's problem}, the result seems to be not optimal for $n\geq2$. Furthermore, it is not trivial but interesting to see how do semilinear damped wave equation and semilinear wave equation affect each other. In this note we only try to understand what happen if we apply directly generalized Kato's lemma.

\medskip

The rest of this note is organized as follows. In Subsection \ref{Section Main result} we show our main results, including local (in time) existence of solutions and blow-up of solutions to \eqref{Nakao's problem}. In Subsection \ref{Section overview} we introduce the overview of our approach, especially, the system of nonlinear ordinary differential inequalities in Lemmas \ref{Lemma system differential} and \ref{Generalzied Kato lemma}. In Sections \ref{Section Lemma System} and \ref{Section Lemma Kato} we prove Lemmas \ref{Lemma system differential} and \ref{Generalzied Kato lemma}, respectively. In Section \ref{Section proof main result} we prove our main results by using the tools developed in the previous sections.

\section{Main results and  overview of our approach}
\setcounter{equation}{0}
\subsection{Main results}\label{Section Main result}
First of all, let us introduce a result of local (in time) existence of Sobolev solutions to the Cauchy problem \eqref{Nakao's problem}. The proof of the following theorem is quite standard (see for example \cite{Strauss1989,Sideris1984,ReissigEbert2018}). Therefore, we will only sketch the proof in Subsection \ref{Subsection proof of local solution}.
\begin{thm}\label{Local existence of solution}
	Let us assume $(u_0,u_1;v_0,v_1)\in \big(H^1(\mb{R}^n)\times L^2(\mb{R}^n)\big)^2$ having compact support, that is, $\text{supp }u_0,u_1,v_0,v_1\subset\{|x|\leq R\}$ with a positive constant $R$.
	\begin{itemize}
		\item If $1<p,q<\infty$ for $n=1,2$, then there exists a positive $T$ and a uniquely determined local (in time) energy solution $(u,v)$ such that
		\begin{align*}
		(u,v)\in\big(\ml{C}\big([0,T],H^1(\mb{R}^n)\big)\cap \ml{C}^1\big([0,T],L^2(\mb{R}^n)\big)\big)^2
		\end{align*}
		with $\text{supp }u(t,\cdot),v(t,\cdot)\subset\{|x|\leq t+R\}$.
		\item If $1<p\leq\frac{n+3}{n-1}$ and $1<q\leq\frac{n}{n-2}$ for $n\geq3$, then there exists a positive $T$ and a uniquely determined local (in time) Sobolev solution $(u,v)$ such that
		\begin{align*}
		(u,v)\in\big(\ml{C}\big([0,T],H^1(\mb{R}^n)\big)\cap \ml{C}^1\big([0,T],L^2(\mb{R}^n)\big)\big)\times\big(\ml{C}\big([0,T],L^{\frac{2(n+1)}{n-1}}(\mb{R}^n)\big)\big)
		\end{align*}
		with $\text{supp }u(t,\cdot),v(t,\cdot)\subset\{|x|\leq t+R\}$.
	\end{itemize}
\end{thm}
	\begin{rem}
	In the cases $n\geq3$ we still can prove local (in time) existence of energy solution
	\begin{align*}
	(u,v)\in\big(\ml{C}\big([0,T],H^1(\mb{R}^n)\big)\cap \ml{C}^1\big([0,T],L^2(\mb{R}^n)\big)\big)^2
	\end{align*}
	with $\text{supp }u(t,\cdot),v(t,\cdot)\subset\{|x|\leq t+R\}$, if $1<p,q\leq\frac{n}{n-2}$. One can see the result in the recent paper \cite{Wakasugi2017}. Nevertheless, we observe that $\frac{n+3}{n-1}>\frac{n}{n-2}$ for all $n\geq3$. For this reason, the consideration of $v\in\ml{C}\big([0,T],L^{\frac{2(n+1)}{n-1}}(\mb{R}^n)\big)$ with $\text{supp }v(t,\cdot)\subset\{|x|\leq t+R\}$ allows us to get a larger admissible range of the exponent $p$.
\end{rem}
Now we state the blow-up results for \eqref{Nakao's problem}. The proof will be shown in Subsection \ref{Subsection Poof blowup}.
\begin{thm}\label{Blow-up result}
	Let $(u_0,u_1;v_0,v_1)\in\big(\ml{C}^2(\mb{R}^n)\times\ml{C}^1(\mb{R}^n)\big)^2$ be nonnegative (but not vanishing) have their support. Assume $(u,v)\in\big(\ml{C}^2([0,T)\times\mb{R}^n)\big)^2$ is the maximal, with respect to time interval, classical solution to the Cauchy problem \eqref{Nakao's problem}. Moreover, we assume that there exists a large time $\tilde{t}_0\in(0,T)$ such that
	\[ \int_{\mb{R}^n}u(\tilde{t}_0,x)\,dx\,\,\,\mbox{ and}\,\,\, \int_{\mb{R}^n}v(\tilde{t}_0,x)\,dx \]	are suitably large (but not infinity). If the exponents satisfy
	\begin{itemize}
		\item $1<p,q<\infty$ for $n=1$,
		\item $1<p,q<\frac{2n}{n-1}$ for $n=2,3$,
		\item $1<p\leq\frac{n+3}{n-1}$, $1<q\leq\frac{n}{n-2}$ for $n\geq4$,
	\end{itemize}
	 and the following condition:
	\begin{align}\label{Condition blow-up}
	\alpha=\max\left\{\frac{q+1}{pq-1};\frac{2+2p^{-1}}{pq-1}\right\}\geq\frac{n-1}{2},
	\end{align}
	then the local (in time) classical solution $(u,v)$ blows up in finite time, that is, $T<\infty$.
\end{thm}
\begin{rem}
Under the assumptions of exponents $p,q$ and the hypothesis for initial data in Theorem \ref{Blow-up result}, according to Theorem \ref{Local existence of solution} we know that there exists a unique local (in time) classical solution to Nakao's problem \eqref{Nakao's problem}.
\end{rem}
\begin{rem}\label{Remark 2.3}
We remark that if we assume suitably large data in Theorem \ref{Blow-up result}, then from the proof of Lemma \ref{Lemma system differential} we may assert that  the time-dependent functions
\begin{align*}
F_1(t)=\int_{\mb{R}^n}u(t,x)\,dx\quad\text{and}\quad F_2(t)=\int_{\mb{R}^n}v(t,x)\,dx
\end{align*}
are strictly increasing functions for $t\geq0$. For this reason, the assumption for suitably large values of \[ F_1(\tilde{t}_0)=\int_{\mb{R}^n}u(\tilde{t}_0,x)\,dx\quad \mbox{and}\quad F_2(\tilde{t}_0)=\int_{\mb{R}^n}v(\tilde{t}_0,x)\,dx \] are trivially satisfied.
\end{rem}
\begin{rem}\label{1D Remark}
	For the one dimensional Nakao's problem, the condition in Theorem \ref{Blow-up result} is always valid and there is no restriction to the exponents of the power nonlinearities.
\end{rem}
\begin{rem}
Actually, the assumption $1<p,q<\frac{2n}{n-1}$ for $n\geq2$ in Theorem \ref{Blow-up result} is a usual assumption when we use the method of proving blow-up for nonlinear ordinary differential inequalities. One may observe this assumption when we prove blow-up of solutions to the single semilinear wave equation (see \cite{Sideris1984,YordanovZhang2006}) and weakly coupled systems of semilinear wave equations (see \cite{DGM,KTW12}). However, due to the symmetry of their models and the proof of the condition $\alpha_{\text{W}}\geq\frac{n-1}{2}$, the condition $1<p,q<\frac{2n}{n-1}$ is trivially valid.
\end{rem}
\begin{rem}
The authors find that the restriction for $F_1(\tilde{t}_0)$ and $F_2(\tilde{t}_0)$ is too strong. We will improve the result in the future.
\end{rem}

\subsection{Overview of our approach}\label{Section overview}
Throughout this note we will consider the time-dependent functions
\begin{align}
F_1(t):&=\int_{\mb{R}^n}u(t,x)dx,\label{F1}\\
F_2(t):&=\int_{\mb{R}^n}v(t,x)dx,\label{F2}
\end{align}
where $(u,v)$ is a classical solution of \eqref{Nakao's problem}.

Since we require compactly supported data, by the property of finite speed of propagation and Theorem \ref{Local existence of solution} it follows that also $(u,v)(t,\cdot)$ is compactly supported with respect to the spatial variables as long as the solution exists with respect to $t$. Thus, if we prove that $F_1(t),F_2(t)$ blow up in finite time, then $(u,v)$ blows up in finite time as well.

In Section \ref{Section Lemma System} we will prove the next lemma by employing the idea of \cite{YordanovZhang2006}.
\begin{lem}\label{Lemma system differential}
	Let us consider the Cauchy problem \eqref{Nakao's problem} with compactly supported initial data with $\text{supp }u_0,u_1,v_0,v_1\subset\{|x|\leq R\}$ for a positive constant $R$. Moreover, we assume
	\begin{align*}
	&\int_{\mb{R}^n}u_j(x)\,dx>0\quad\text{and}\quad \int_{\mb{R}^n}v_j(x)\,dx>0\quad\text{for}\quad j=0,1.
	\end{align*}
	Moreover, we assume  $1<p,q<\infty$ for $n=1$, and $1<p,q<\frac{2n}{n-1}$ for $n\geq 2$.
	Then, the functions $F_1(t)$, $F_2(t)$ satisfy the following system of second-order nonlinear differential inequalities:
	\begin{align*}
	F_1(t)&\geq C_3(t+R)^{1+\frac{2-p}{2}(n-1)},\\
	\frac{dF_1}{dt}(t)+F_1(t)&\geq 4C_3(t+R)^{1+\frac{2-p}{2}(n-1)},\\
	\frac{d^2F_1}{dt^2}(t)+\frac{dF_1}{dt}(t)&\geq (t+R)^{-n(p-1)}(F_2(t))^p,\\
	F_2(t)&\geq (t+R),\\
	\frac{d^2F_2}{dt^2}(t)&\geq e^{-\frac{3-\sqrt{5}}{2}qt}(t+R)^{-n(q-1)}(F_1(t))^q,
	\end{align*}
	with $C_3>0$ defined in \eqref{C_3}, respectively, for large time $t\geq T_0$.
\end{lem}
\begin{rem}
To investigate the lower bounds estimate of $F_1(t)$ in Lemma \ref{Lemma system differential}, we need to assume $p<\frac{2n}{n-1}$ when $n\geq2$. Similarly, to get the lower bounds estimate of $F_2(t)$ in Lemma \ref{Lemma system differential}, we need to assume $q<\frac{2n}{n-1}$ when $n\geq2$.
\end{rem}

Then, in Section \ref{Section Lemma Kato} we may prove blow-up of the functions $F_1(t),F_2(t)$ in finite time by using the following generalized Kato's lemma.
\begin{lem}\label{Generalzied Kato lemma}
Let $p,q>1$, $\alpha_1,\beta_1>0$, $\alpha_2,\beta_2,\beta_3\geq0$. Moreover, we assume
\begin{align}
\beta_2+\alpha_2q&\leq\beta_1(pq-1)+2(q+1),\label{Cond 1}\\
\alpha_2+\beta_2p&\leq\alpha_1(pq-1)+2(p+1).\label{Cond 2}
\end{align}
 Suppose that $F_1=F_1(t)$ and $F_2=F_2(t)$ belong to $\ml{C}^2[0,T]$ and satisfy the following system of second-order nonlinear differential inequalities for $t\geq T_0$:
\begin{align}
F_1(t)&\geq k_0(t+R)^{\alpha_1},\label{1EQ 2}\\
\frac{dF_1}{dt}(t)+F_1(t)&\geq k_1(t+R)^{\alpha_1},\label{1EQ 0}\\
\frac{d^2F_1}{dt^2}(t)+\frac{dF_1}{dt}(t)&\geq k_2(t+R)^{-\alpha_2}(F_2(t))^p,\label{1EQ 1}\\
F_2(t)&\geq k_3(t+R)^{\beta_1},\label{1EQ 4}\\
\frac{d^2F_2}{dt^2}(t)&\geq k_4e^{-\beta_3t}(t+R)^{-\beta_2}(F_1(t))^q,\label{1EQ 3}
\end{align}
where all $k_0,k_1,k_2,k_3,k_4$ are positive constants.
Moreover, we assume that there exists a large $\widetilde{T}_0$ such that $\widetilde{T}_0\in(T_0,T)$ such that $F_1\big(\widetilde{T}_0\big)$ and $F_2\big(\widetilde{T}_0\big)$ are suitably large (but not infinity). Then, the functions $F_1(t),F_2(t)$ blow up in finite time, that is, $T<\infty$.
\end{lem}
\section{Proof of Lemma \ref{Lemma system differential}}\label{Section Lemma System}
\setcounter{equation}{0}
\subsection{Derivation of a system of nonlinear ordinary differential inequalities} \label{Secdifferentialinequalities}
Let $(u,v)$ be the local (in time) classical solutions to the Cauchy problem \eqref{Nakao's problem}. We define the functions $F_1=F_1(t)$ and $F_2=F_2(t)$ as in \eqref{F1} and \eqref{F2}, respectively. From the property of finite speed of propagation and our assumption on the initial data, we have
\begin{align*}
\text{supp }u(t,\cdot),v(t,\cdot)\subset\{|x|\leq t+R\}\quad\text{ for any }t>0.
\end{align*}
By applying the divergence theorem and the compact support property of solutions, we have
\begin{align*}
\frac{d^2F_1}{dt^2}(t)+\frac{dF_1}{dt}(t)&=\int_{\mb{R}^n}\big(u_{tt}(t,x)+u_t(t,x)\big)\,dx=\int_{\mb{R}^n}\big(\Delta u(t,x)+|v(t,x)|^p\big)\,dx=\int_{\mb{R}^n}|v(t,x)|^p\,dx,\\
\frac{d^2F_2}{dt^2}(t)&=\int_{\mb{R}^n}v_{tt}(t,x)\,dx=\int_{\mb{R}^n}\big(\Delta v(t,x)+|u(t,x)|^q\big)\,dx=\int_{\mb{R}^n}|u(t,x)|^q\,dx.
\end{align*}
Using H\"older's inequality leads to
\begin{align*}
\Big|\int_{\mb{R}^n}v(t,x)\,dx\Big|^p&\leq\Big(\int_{|x|\leq t+R}|v(t,x)|^p\,dx\Big)\Big(\int_{|x|\leq t+R}1\,dx\Big)^{p-1},\\
\Big|\int_{\mb{R}^n}u(t,x)\,dx\Big|^q&\leq\Big(\int_{|x|\leq t+R}|u(t,x)|^q\,dx\Big)\Big(\int_{|x|\leq t+R}1\,dx\Big)^{q-1}.
\end{align*}
In other words, we obtain
\begin{align*}
\frac{d^2 F_1}{dt^2}(t)+\frac{dF_1}{dt^2}(t)&\geq C (t+R)^{-n(p-1)}|F_2(t)|^p,\\
\frac{d^2 F_2}{dt^2}(t)&\geq C (t+R)^{-n(q-1)}|F_1(t)|^q.
\end{align*}
Here and in the following we use $C$ as a universal positive constant.
In order to get  lower bounds for the functions $F_1=F_1(t)$ and $F_2=F_2(t)$, we now introduce the functions $\psi_1=\psi_1(t,x)$ and $\psi_2=\psi_2(t,x)$ as follows:
\begin{align*}
\psi_1(t,x):=e^{-\frac{\sqrt{5}-1}{2}t}\phi(x)\quad\text{and}\quad \psi_2(t,x):=e^{-t}\phi(x),
\end{align*}
where
\begin{align*}
\phi(x):=\int_{\mb{S}^{n-1}}e^{x\cdot\omega}d\sigma_{\omega}\sim |x|^{-\frac{n-1}{2}}e^{|x|}\quad\text{as}\quad|x|\rightarrow\infty,
\end{align*}
and $\mb{S}^{n-1}$ is the $n-1$ dimensional sphere. By the compactness of the unit sphere $\mb{S}^{n-1}$ we know that the function $\phi \in\ml{C}^{\infty}(\mb{R}^n)$ satisfies
\begin{align*}
\phi(x)-\Delta\phi(x)=\phi(x)-\int_{\mb{S}^{n-1}}\sum\limits_{k=1}^n\omega_k^2e^{x\cdot\omega}d\sigma_{\omega}=0.
\end{align*}
Actually, the fact that $\phi(x)>0$ for $x\in\mb{R}^n$ comes from the monotonicity of the exponential function. Moreover, from direct calculations, we may assert that
\begin{align}
\partial_t^2\psi_1(t,x)-\Delta\psi_1(t,x)-\partial_t\psi_1(t,x)&=0,\label{Psi 1}\\
\partial_t^2\psi_2(t,x)-\Delta\psi_2(t,x)&=0.\label{Psi 2}
\end{align}
We apply reverse H\"older's inequality to obtain
\begin{align*}
\frac{d^2 F_1}{dt^2}(t)+\frac{d F_1}{dt}(t)\geq&\Big(\int_{\mb{R}^n}v(t,x)\psi_2(t,x)\,dx\Big)^p\Big(\int_{|x|\leq t+R}\psi_2(t,x)^{\frac{p}{p-1}}\,dx\Big)^{-(p-1)}=:J_1(t)\, J_2(t),\\
\frac{d^2F_2}{dt^2}(t)\geq&\Big(\int_{\mb{R}^n}u(t,x)\psi_1(t,x)\,dx\Big)^q\Big(\int_{|x|\leq t+R}\psi_1(t,x)^{\frac{q}{q-1}}\,dx\Big)^{-(q-1)}=:J_3(t)\, J_4(t).
\end{align*}
In the following steps we will estimate the time-dependent functions $J_1=J_1(t)$, $J_3=J_3(t)$ and $J_2=J_2(t)$, $J_4=J_4(t)$, respectively.
\subsection{Estimate of the time-dependent function $J_1$} \label{SecestimateJ1}
\noindent In this step we follow \cite{YordanovZhang2006}. Multiplying the equation \eqref{Nakao's problem}$_2$ by the function $\psi_2$ and integrating it over $[0,\eta]\times\mb{R}^n$, we obtain
\begin{align*}
0\leq&\int_0^{\eta}\int_{\mb{R}^n}|u(t,x)|^q\psi_2(t,x)\,dx\,dt\\
&=\int_0^{\eta}\int_{\mb{R}^n}\big(v_{tt}(t,x)-\Delta v(t,x)\big)\psi_2(t,x)\,dx\,dt\\
&=\int_{0}^{\eta}\int_{\mb{R}^n}v(t,x)\big(\partial_t^2\psi_{2}(t,x)-\Delta\psi_2(t,x)\big)\,dx\,dt\\
&\quad\,\,+\int_{\mb{R}^n}(v_t(t,x)\psi_2(t,x)-v(t,x)\partial_t\psi_2(t,x))\,dx\Big|_{t=0}^{t=\eta}\\
&=\int_{\mb{R}^n}\big(v_{\eta}(\eta,x)\psi_2(\eta,x)-v(\eta,x)\partial_{\eta}\psi_2(\eta,x)\big)\,dx-\int_{\mb{R}^n}\big(v_0(x)+v_1(x)\big)\phi(x))\,dx,
\end{align*}
where we used $\psi_2(t,x)\geq0$ and the equation \eqref{Psi 2}.
We now define the time-dependent function
\begin{align*}
F_3(\eta):=\int_{\mb{R}^n}v(\eta,x)\psi_2(\eta,x)\,dx.
\end{align*}
Thus,
\begin{align}\label{Supp 5}
\frac{dF_3}{d\eta}(\eta)+2F_3(\eta)\geq\int_{\mb{R}^n}\big(v_0(x)+v_1(x)\big)\phi(x)\,dx,
\end{align}
where we used the relation $\partial_{\eta}\psi_2(\eta,x)=-\psi_2(\eta,x)$.
The assumption on the initial data implies
\begin{align*}
\int_{\mb{R}^n}\big(v_0(x)+v_1(x)\big)\phi(x)\,dx>0\quad\text{and}\quad\int_{\mb{R}^n}v_0(x)\phi(x)\,dx>0.
\end{align*}
Then, multiplying \eqref{Supp 5} by $e^{2\eta}$ and integrating it over $[0,t]$, we obtain
\begin{align*}
F_3(t)\geq e^{-2t}\int_{\mb{R}^n}v_0(x)\phi(x)\,dx+\frac{1}{2}\big(1-e^{-2t}\big)\int_{\mb{R}^n}\big(v_0(x)+v_1(x)\big)\phi(x)\,dx\geq C_0,
\end{align*}
with a positive constant $C_0$. We immediately conclude $J_1(t)\geq C_0^p$.
\subsection{Estimate of the time-dependent function $J_3$} \label{SecestimateJ3}
\noindent Noting that $\psi_1(t,x)\geq0$ we multiply the equation \eqref{Nakao's problem}$_1$ by the function $\psi_1$ and integrate it over $[0,\eta]\times\mb{R}^n$ to derive
\begin{align*}
0\leq&\int_0^{\eta}\int_{\mb{R}^n}|v(t,x)|^p\psi_1(t,x)\,dx\,dt\\
&=\int_0^{\eta}\int_{\mb{R}^n}\big(u_{tt}(t,x)-\Delta u(t,x)+u_t(t,x)\big)\psi_1(t,x)\,dx\,dt\\
&=\int_{0}^{\eta}\int_{\mb{R}^n}u(t,x)\big(\partial_t^2\psi_{1}(t,x)-\Delta\psi_1(t,x)-\partial_t\psi_1(t,x)\big)\,dx\,dt\\
&\quad\,\,+\int_{\mb{R}^n}\big(u_t(t,x)\psi_1(t,x)-u(t,x)\partial_t\psi_1(t,x)+u(t,x)\psi_1(t,x)\big)\,dx\Big|_{t=0}^{t=\eta}.
\end{align*}
Here we define the time-dependent function
\begin{equation*}
F_4(\eta):=\int_{\mb{R}^n}u(\eta,x)\psi_1(\eta,x)\,dx.
\end{equation*}
By using the equation \eqref{Psi 1} and the relation
\begin{align*}
& u_t(t,x)\psi_1(t,x)-u(t,x)\partial_t\psi_1(t,x)+u(t,x)\psi_1(t,x)\\ & \qquad =\partial_t (u(t,x)\psi_1(t,x))+\sqrt{5}u(t,x)\psi_1(t,x),
\end{align*}
we have
\begin{align}\label{Supp 01}
\frac{d}{d\eta}F_4(\eta)+\sqrt{5}F_4(\eta)\geq\int_{\mb{R}^n}\Big(u_1(x)+\frac{1+\sqrt{5}}{2}u_0(x)\Big)\phi(x)\,dx.
\end{align}
Applying Gronwall's inequality to \eqref{Supp 01} implies
\begin{align*}
F_4(t)\geq e^{-\sqrt{5}t}\int_{\mb{R}^n}u_0(x)\phi(x)\,dx+\frac{1}{\sqrt{5}}\big(1-e^{-\sqrt{5}t}\big)\int_{\mb{R}^n}\Big(u_1(x)+\frac{1+\sqrt{5}}{2}u_0(x)\Big)\phi(x)\,dx\geq C_1,
\end{align*}
with a positive constant $C_1$, where we used again the assumption for the initial data. We can immediately conclude $J_3(t)\geq C_1^{q}$.
\subsection{Estimate of the time-dependent functions $J_2$ and $J_4$} \label{SecestimateJ2J4}
\noindent For one thing, we may use the asymptotic behavior of $\phi$ to get
\begin{equation*}
\int_{|x|\leq t+R}\psi_2(t,x)^{\frac{r}{r-1}}\,dx\leq C_{2} (t+R)^{n-1-\frac{(n-1)p'}{2}},
\end{equation*}
with a positive constant $C_{2}$. \\
For another, by using asymptotic behavior of $\phi$ again we obtain
\begin{align*}
\int_{|x|\leq t+R}\psi_1(t,x)^{\frac{r}{r-1}}\,dx&\leq \widetilde{C}_0e^{-\frac{(\sqrt{5}-1)q'}{2}t}\int_{|x|\leq t+R}|x|^{-\frac{(n-1)q'}{2}}e^{|x|q'}dx\\
&\leq\widetilde{C}_0e^{-\frac{(\sqrt{5}-1)q'}{2}t}\int_{0}^{t+R}r^{-\frac{(n-1)q'}{2}+n-1}e^{rq'}dr\\
&\leq\frac{1}{q'}\widetilde{C}_0e^{-\frac{(\sqrt{5}-1)q'}{2}t}\Big(e^{(t+R)q'}(t+R)^{n-1-\frac{(n-1)q'}{2}}-1\\
&\qquad\qquad\qquad\qquad\quad-\Big(n-1-\frac{(n-1)q'}{2}\Big)\int_0^{t+R}e^{rq'}r^{n-2-\frac{(n-1)q'}{2}}dr\Big)\\
&\leq\frac{1}{q'}\widetilde{C}_0e^{-\frac{(\sqrt{5}-1)q'}{2}t+(t+R)q'}(t+R)^{n-1-\frac{(n-1)q'}{2}}\\
&=:\widetilde{C}_2e^{\frac{3-\sqrt{5}}{2}q't}(t+R)^{n-1-\frac{(n-1)q'}{2}},
\end{align*}
with a positive constant $\widetilde{C}_2$.

 It shows that
\begin{align*}
J_2(t)&\geq C_{2}^{-(p-1)} (t+R)^{\frac{2-p}{2}(n-1)},\\
 J_4(t)&\geq \widetilde{C}_{2}^{-(q-1)}e^{-\frac{3-\sqrt{5}}{2}qt} (t+R)^{\frac{2-q}{2}(n-1)}.
\end{align*}
Summarizing the derived estimates from the above subsections concludes
\begin{align}
\frac{d^2 F_1}{dt^2}(t)+\frac{dF_1}{dt}(t)&\geq C_0^pC_{2}^{-(p-1)}(t+R)^{\frac{2-p}{2}(n-1)},\label{Lower bound F1 01}\\
\frac{d^2F_2}{dt^2}(t)&\geq C_1^q \widetilde{C}_{2}^{-(q-1)}e^{-\frac{3-\sqrt{5}}{2}qt}(t+R)^{\frac{2-q}{2}(n-1)}\label{Lower bound F2 01}.
\end{align}
\subsection{Lower bound for $F_1$} \label{SeclowerboundF1}
Let us derive a lower bound for $F_1=F_1(t)$ by using \eqref{Lower bound F1 01}. Integrating \eqref{Lower bound F1 01} over $[0,t]$ gives
\begin{equation}\label{Supp 2}
\begin{aligned}
\frac{d F_1}{dt}(t)+F_1(t)&\geq \frac{d F_1}{dt}(0)+F_1(0)+\frac{C_0^pC_{2}^{-(p-1)}}{1+\frac{2-p}{2}(n-1)}\Big((t+R)^{1+\frac{2-p}{2}(n-1)}-R^{1+\frac{2-p}{2}(n-1)}\Big)\\
&\geq \frac{d F_1}{dt}(0)+F_1(0)+\frac{C_0^pC_{2}^{-(p-1)}}{2+(2-p)(n-1)}(t+R)^{1+\frac{2-p}{2}(n-1)},
\end{aligned}
\end{equation}
for large time $t\geq t_0$, where we used our assumption
\begin{align*}
\mbox{from}\,\,\,1<p<\frac{2n}{n-1}\quad\text{if}\quad n\geq2\quad\mbox{it follows}\quad 1+\frac{2-p}{2}(n-1)>0\quad\text{for all}\quad n\geq1.
\end{align*} Then, we use Gronwall's inequality to \eqref{Supp 2} to find
\begin{align*}
F_1(t)\geq &\big(1-e^{t_0-t}\big)\Big(\frac{dF_1}{dt}(0)+F_1(0)\Big)+e^{t_0-t}F_1(t_0)\\
&+\frac{C_0^pC_{2}^{-(p-1)}}{2+(2-p)(n-1)}e^{-t}\int_{t_0}^te^{\tau}(\tau+R)^{1+\frac{2-p}{2}(n-1)}\,d\tau.
\end{align*}
The use of integration by parts implies
\begin{align*}
&e^{-t}\int_{t_0}^te^{\tau}(\tau+R)^{1+\frac{2-p}{2}(n-1)}d\tau\\
&=(t+R)^{1+\frac{2-p}{2}(n-1)}-e^{t_0-t}(t_0+R)^{1+\frac{2-p}{2}(n-1)}-\Big(1+\frac{2-p}{2}(n-1)\Big)(t+R)^{\frac{2-p}{2}(n-1)}\\
&\quad\,\,+\Big(1+\frac{2-p}{2}(n-1)\Big) e^{t_0-t}(t_0+R)^{\frac{2-p}{2}(n-1)}\\
&\quad\,\,+\frac{2-p}{2}(n-1)\Big(1+\frac{2-p}{2}(n-1)\Big)e^{-t}\int_{t_0}^te^{\tau}(\tau+R)^{\frac{2-p}{2}(n-1)-1}\,d\tau,\\
&\geq\frac{1}{4} (t+R)^{1+\frac{2-p}{2}(n-1)},
\end{align*}
where we applied the following facts:
\begin{align*}
e^{t_0-t}(t_0+R)^{1+\frac{2-p}{2}(n-1)}&\leq \frac{1}{4}(t+R)^{1+\frac{2-p}{2}(n-1)},\\
\Big(1+\frac{2-p}{2}(n-1)\Big)(t+R)^{\frac{2-p}{2}(n-1)}&\leq \frac{1}{4}(t+R)^{1+\frac{2-p}{2}(n-1)},
\end{align*}
and
\begin{align*}
&\Big|\frac{2-p}{2}(n-1)\Big(1+\frac{2-p}{2}(n-1)\Big)e^{t_0-t}\int_{t_0}^te^{\tau}(\tau+R)^{\frac{2-p}{2}(n-1)-1}\,d\tau\Big|\\
&\leq\Big|\frac{2-p}{2}(n-1)\Big|\Big(1+\frac{2-p}{2}(n-1)\Big)e^{t_0}\Big|\int_{t_0}^t(\tau+R)^{\frac{2-p}{2}(n-1)-1}\,d\tau\Big|\\
&\leq \frac{1}{4}(t+R)^{1+\frac{2-p}{2}(n-1)},
\end{align*}
for $t\geq t_1>t_0$, where $t_1$ is sufficiently large. So, we immediately get for $t\geq t_1$ the estimate to below
\begin{align*}
F_1(t)\geq C_3(t+R)^{1+\frac{2-p}{2}(n-1)},
\end{align*}
where
\begin{align}\label{C_3}
C_3:=\frac{C_0^pC_{2}^{-(p-1)}}{4(2+(2-p)(n-1))}.
\end{align}
\subsection{Lower bound for $F_2$} \label{SeclowerboundF2}
To get a lower bound for $F_2=F_2(t)$, we only need to integrate twice with respect to  $t$ the inequality \eqref{Lower bound F2 01}. In this way we obtain
\begin{equation}\label{Supp 9}
\begin{aligned}
F_2(t)\geq&F_2(t_2)+(t-t_2)\frac{dF_2}{dt}(0)\\
&+\frac{C_1^q C_{2}^{-(q-1)}e^{-\frac{3-\sqrt{5}}{2}qt}}{2\big(1+\frac{2-q}{2}(n-1)\big)\big(2+\frac{2-q}{2}(n-1)\big)}\big((t+R)^{2+\frac{2-q}{2}(n-1)}-(t_2+R)^{2+\frac{2-q}{2}(n-1)}\big),
\end{aligned}
\end{equation}
for $t\geq t_2$. Under our assumption $1<q<\frac{2n}{n-1}$ if $n\geq2$ the following estimate holds for $t\geq t_3>t_2$:
\begin{align*}
F_2(t)\geq (t+R).
\end{align*}
Here $t_3$ is a sufficiently large positive constant.
In conclusion, taking $T_0=\max\{t_1;t_3\}$ we derived all desired estimates. The proof is complete.
\section{Proof of Lemma \ref{Generalzied Kato lemma}}\label{Section Lemma Kato}
\setcounter{equation}{0}
Let us describe some properties for the functions $F_1=F_1(t)$ and $F_2=F_2(t)$. From \eqref{1EQ 3}, \eqref{1EQ 4} and \eqref{1EQ 2} we obtain $\frac{d^2F_2}{dt^2}(t)>0$ and $F_2(t)>0$ for all $t\geq T_0$. Thus, $F_2$ is a convex function for $t\geq T_0$, which implies that there exists $T_1\geq T_0$ such that $\frac{dF_2}{dt}(t)>0$ for $t\geq T_1$. Similarly, from \eqref{1EQ 0} and \eqref{1EQ 2} we know $\frac{dF_1}{dt}(t)+F_1(t)>0$ and $F_1(t)>0$ for all $t\geq T_0$. Let us apply Gronwall's inequality to \eqref{1EQ 1} together with \eqref{1EQ 4} for  we also have
 \begin{align*}
 \frac{dF_1}{dt}(t)&\geq e^{T_0-t}\frac{dF_1}{dt}(T_0)+k_2k_3^p\inf\limits_{\tau\in[T_0,t]}(\tau+R)^{-\alpha_2+\beta_1p}\,(1-e^{T_0-t})>0,
 \end{align*}
 for $t\geq T_2>T_0$. All in all, we have
 \begin{align*}
 \frac{dF_1}{dt}(t)>0\quad\text{and}\quad \frac{dF_2}{dt}(t)>0,
 \end{align*}
 for large time $t\geq T_3:=\max\{T_1,T_2\}$.\\
To prove Lemma \ref{Generalzied Kato lemma}, we discuss the finite time blow-up of the time-dependent function $F_2$ in Subsection \ref{Subsection F2}, and the finite time blow-up of the time-dependent function $F_1$ in Subsection \ref{Subsection F1}.
\subsection{Finite time blow-up of $F_2=F_2(t)$}\label{Subsection F2}
 First of all, multiplying \eqref{1EQ 1} by $\frac{dF_2}{dt}(t)$ we obtain
\begin{align}\label{EQ 5}
\frac{d}{dt}\Big(\frac{dF_1}{dt}(t)+F_1(t)\Big)\,\dfrac{dF_2}{dt}(t)\geq\frac{k_2}{p+1}(t+R)^{-\alpha_2}\frac{d}{dt}\big((F_2(t))^{p+1}\big).
\end{align}
Applying integration by parts leads to
\begin{align*}
\int_{T_3}^t\frac{d}{d\tau}\Big(\frac{dF_1}{d\tau}(\tau)+F_1(\tau)\Big)\,\dfrac{dF_2}{d\tau}(\tau)\,d\tau&=\Big(\frac{dF_1}{dt}(t)+F_1(t)\Big)\dfrac{dF_2}{dt}(t)-\Big(\frac{dF_1}{dt}(T_3)+F_1(T_3)\Big)\dfrac{dF_2}{dt}(T_3)\\
&\quad\,\,-\int_{T_3}^t\Big(\frac{dF_1}{d\tau}(\tau)+F_1(\tau)\Big)\cdot\dfrac{d^2F_2}{d\tau^2}(\tau)\,d\tau\\
&\leq \Big(\frac{dF_1}{dt}(t)+F_1(t)\Big)\dfrac{dF_2}{dt}(t).
\end{align*}
Next, we integrate \eqref{EQ 5} over $[T_3,t]$ to have
\begin{equation}\label{EQ 11}
\begin{aligned}
\Big(\frac{dF_1}{dt}(t)+F_1(t)\Big)\frac{dF_2}{dt}(t)&\geq\frac{k_2}{p+1}\inf_{\tau\in[T_3,t]}(\tau+R)^{-\alpha_2}\cdot\big((F_2(t))^{p+1}-(F_2(T_3))^{p+1}\big)\\
&\geq\frac{k_2}{2(p+1)}(t+R)^{-\alpha_2}(F_2(t))^{p+1},
\end{aligned}
\end{equation}
due to the monotonically increasing property of $F_2$, where $t\geq T_4> T_3$ with a sufficiently large $T_4$.\\
Then, we multiply the above inequality by $e^t\frac{dF_2}{dt}$ to get
\begin{equation}\label{EQ 16}
\begin{aligned}
\frac{d}{dt}\big(e^tF_1(t)\big)\,\Big(\frac{dF_2}{dt}(t)\Big)^2&=e^t\Big(\frac{dF_1}{dt}(t)+F_1(t)\Big)\,\Big(\frac{dF_2}{dt}(t)\Big)^2\\
&\geq \frac{k_2}{2(p+1)(p+2)}e^t(t+R)^{-\alpha_2}\frac{d}{dt}\big((F_2(t))^{p+2}\big).
\end{aligned}
\end{equation}
Similarly, we know from integration by parts
\begin{align*}
\int_{T_4}^t\frac{d}{d\tau}\big(e^\tau F_1(\tau)\big)\,\Big(\frac{dF_2}{d\tau}(\tau)\Big)^2d\tau&=e^tF_1(t)\Big(\frac{dF_2}{dt}(t)\Big)^2-e^{T_4}F_1(T_4)\Big(\frac{dF_2}{dt}(T_4)\Big)^2\\
&\quad\,\,-2\int_{T_4}^te^\tau F_1(\tau)\frac{dF_2}{d\tau}(\tau)\,\frac{d^2F_2}{d\tau^2}(\tau)d\tau\\
&\leq e^tF_1(t)\Big(\frac{dF_2}{dt}(t)\Big)^2.
\end{align*}
Again by using the monotonically increasing behavior of $F_2$ for $t\geq T_5>T_4$, integrating \eqref{EQ 16} over $[T_4,t]$ again shows
\begin{align*}
e^tF_1(t)\Big(\frac{dF_2}{dt}(t)\Big)^2&\geq\frac{k_2}{2(p+1)(p+2)}\inf_{\tau\in[T_4,t]}e^{\tau}\,\inf_{\tau\in[T_4,t]}(\tau+R)^{-\alpha_2}\cdot\big((F_2(t))^{p+2}-(F_2(T_4))^{p+2}\big)\\
&\geq\frac{k_2}{4(p+1)(p+2)}(t+R)^{-\alpha_2}(F_2(t))^{p+2}.
\end{align*}
Thus, we derived the following lower bound estimate for $F_1(t)$:
\begin{align}\label{EQ 7}
F_1(t)\geq \frac{k_2}{4(p+1)(p+2)}e^{-t}(t+R)^{-\alpha_2}(F_2(t))^{p+2}\Big(\frac{dF_2}{dt}(t)\Big)^{-2}.
\end{align}
Plugging \eqref{EQ 7} into \eqref{1EQ 3} implies
\begin{align*}
\frac{1}{2q+1}\frac{d}{dt}\Big(\Big(\frac{dF_2}{dt}(t)\Big)^{2q+1}\Big)&=\frac{d^2F_2}{dt^2}(t)\cdot\Big(\frac{dF_2}{dt}(t)\Big)^{2q}\\
&\geq \frac{k_2^qk_4}{(4(p+1)(p+2))^q}e^{-(q+\beta_3)t}(t+R)^{-\beta_2-\alpha_2q}(F_2(t))^{(p+2)q}.
\end{align*}
Multiplying the above inequality by $\frac{dF_2}{dt}$ once more, we have
\begin{align*}
\frac{dF_2}{dt}(t)\,\frac{d}{dt}\Big(\Big(\frac{dF_2}{dt}(t)\Big)^{2q+1}\Big)\geq& 2k_5e^{-(q+\beta_3)t}(t+R)^{-\beta_2-\alpha_2q}\frac{d}{dt}\big((F_2(t))^{(p+2)q+1}\big),
\end{align*}
where
\begin{align*}
k_5:=\frac{(2q+1)k_2^qk_4}{2((p+2)q+1)(4(p+1)(p+2))^q}.
\end{align*}
Therefore, combining
\begin{align*}
&\int_{T_5}^t\frac{dF_2}{d\tau}(\tau)\,\frac{d}{d\tau}\Big(\Big(\frac{dF_2}{d\tau}(\tau)\Big)^{2q+1}\Big)\,d\tau\\
&=\Big(\frac{dF_2}{dt}(t)\Big)^{2(q+1)}-\Big(\frac{dF_2}{dt}(T_5)\Big)^{2(q+1)}-\int_{T_5}^t\frac{d^2F_2}{d\tau^2}(\tau)\Big(\frac{dF_2}{d\tau}(\tau)\Big)^{2q+1}\,d\tau\\
&\leq \Big(\frac{dF_2}{dt}(t)\Big)^{2(q+1)},
\end{align*}
 with
\begin{align*}
&\int_{T_5}^t(\tau+R)^{-\beta_2-\alpha_2q}\frac{d}{dt}\big((F_2(\tau))^{(p+2)q+1}\big)\,d\tau\\
&\geq\inf\limits_{\tau\in[T_5,t]}\big(e^{-q\tau}(\tau+R)^{-\beta_2-\alpha_2q}\big)\,\big((F_2(t))^{(p+2)q+1}-(F_2(T_5))^{(p+2)q+1}\big)\\
&\geq\frac{1}{2}e^{-qt}(t+R)^{-\beta_2-\alpha_2q}(F_2(t))^{(p+2)q+1}
\end{align*}
for $t\geq T_6>T_5$, the following estimate holds:
\begin{align}\label{EQ 8}
\frac{dF_2}{dt}(t)\geq k_5^{\frac{1}{2(q+1)}}e^{-\frac{q+\beta_3}{2(q+1)}t}(t+R)^{-\frac{\beta_2+\alpha_2q}{2(q+1)}}(F_2(t))^{\frac{(p+2)q+1}{2(q+1)}}.
\end{align}
Here $T_6$ is a large positive constant.

To prove our desired statements, we shall distinguish between two cases.
\subsubsection{The condition $\beta_2+\alpha_2q<\beta_1(pq-1)+2(q+1)$ holds.}

\noindent The estimate \eqref{1EQ 4} shows that
\begin{align}\label{Supp 6}
(F_2(t))^{\epsilon}\geq k_3^{\epsilon}(t+R)^{\beta_1\epsilon},
\end{align}
for a constant $\epsilon>0$ to be determined later. Then, plugging \eqref{Supp 6} into \eqref{EQ 8} yields
\begin{align}\label{Supp 7}
\frac{dF_2}{dt}(t)\geq k_3^{\epsilon}k_5^{\frac{1}{2(q+1)}}e^{-\frac{q+\beta_3}{2(q+1)}t}(t+R)^{-\frac{\beta_2+\alpha_2q}{2(q+1)}+\beta_1\epsilon}(F_2(t))^{\frac{(p+2)q+1}{2(q+1)}-\epsilon}.
\end{align}
Obviously, our assumption
\begin{align*}
\beta_2+\alpha_2q<\beta_1(pq-1)+2(q+1)
\end{align*}
can be rewritten by
\begin{align*}
\left(\frac{\beta_2+\alpha_2q}{2(q+1)}-1\right)\frac{1}{\beta_1}<\epsilon<\frac{(p+2)q+1}{2(q+1)}-1,
\end{align*}
 which means that
\begin{align*}
-\frac{\beta_2+\alpha_2q}{2(q+1)}+\beta_1\epsilon>-1
\quad\text{and}\quad\frac{(p+2)q+1}{2(q+1)}-\epsilon>1.
\end{align*}
Let us consider the auxiliary initial value problem
\begin{align}\label{Y(t)}
\frac{dY}{dt}(t)=\kappa e^{-\nu t}(t+R)^{-\alpha}(Y(t))^{\beta},\quad Y(T_6) = F_2(T_6),
\end{align}
with $\kappa>0$ and $\nu>0$. The solution of \eqref{Y(t)} is given by
\begin{align*}
Y(t)=\Big((Y(T_6))^{1-\beta}-\kappa(\beta-1)\int_{T_6}^te^{-\nu\tau}(\tau+R)^{-\alpha}d\tau\Big)^{-\frac{1}{\beta-1}}.
\end{align*}
It is clear that if $-\alpha>-1$ and $\beta>1$. Then the solution of \eqref{Y(t)} blows up when  $Y(T_6)$ is large. According to the Petrovitsch theorem, we conclude that $F_2=F_2(t)$ blows up in finite time.
\subsubsection{The condition $\beta_2+\alpha_2q=\beta_1(pq-1)+2(q+1)$ holds.} \label{Sec1000}
\noindent In this case we choose a positive constant $\delta$ such that $\delta\in\big(0,\big(\frac{\beta_2+\alpha_2q}{2(q+1)}-1\big)\frac{1}{\beta_1}\big)$. Then we multiply \eqref{EQ 8} by $(F_2(t))^{-1-\delta}$ to get immediately
\begin{equation}\label{Supp 8}
\begin{aligned}
-\frac{1}{\delta}\frac{d}{dt}\big((F_2(t))^{-\delta}\big)&\geq k_5^{\frac{1}{2(q+1)}}(t+R)^{-\frac{\beta_2+\alpha_2q}{2(q+1)}}e^{-\frac{q+\beta_3}{2(q+1)}t}
(F_2(t))^{-\delta+(\frac{\beta_2+\alpha_2q}{2(q+1)}-1)\frac{1}{\beta_1}}\\
&\geq k_3^{-\delta+(\frac{\beta_2+\alpha_2q}{2(q+1)}-1)\frac{1}{\beta_1}}k_5^{\frac{1}{2(q+1)}}e^{-\frac{q+\beta_3}{2(q+1)}t}(t+R)^{-1-\delta\beta_1},
\end{aligned}
\end{equation}
 where we used our estimate \eqref{1EQ 4} and our basic assumption in this case which reads as follows:
\begin{align*}
\Big(\frac{\beta_2+\alpha_2q}{2(q+1)}-1\Big)\frac{1}{\beta_1}=\frac{(p+2)q+1}{2(q+1)}-1.
\end{align*}
Thus, we integrate \eqref{Supp 8} over $[T_6,t]$ to get
\begin{align*}
(F_2(T_6))^{-\delta}-(F_2(t))^{-\delta}\geq&\frac{1}{\beta_1} k_3^{-\delta+\big(\frac{\beta_2+\alpha_2q}{2(q+1)}-1\big)\frac{1}{\beta_1}}k_5^{\frac{1}{2(q+1)}}e^{-\frac{q+\beta_3}{2(q+1)}t}\big((T_6+R)^{-\delta\beta_1}-(t+R)^{-\delta\beta_1}\big)\\
&=k_6e^{-\frac{q}{2(q+1)}t}\big((T_6+R)^{-\delta\beta_1}-(t+R)^{-\delta\beta_1}\big),
\end{align*}
where the constant
\begin{align*}
k_6:=\frac{1}{\beta_1} k_3^{-\delta+\big(\frac{\beta_2+\alpha_2q}{2(q+1)}-1\big)\frac{1}{\beta_1}}k_5^{\frac{1}{2(q+1)}}
\end{align*}
is of course independent of $T_6$. In other words, we have
\begin{align*}
F_2(t)\geq\Big((F_2(T_6))^{-\delta}-k_6e^{-\frac{q+\beta_3}{2(q+1)}t}\big((T_6+R)^{-\delta\beta_1}-(t+R)^{-\delta\beta_1}\big)\Big)^{-\frac{1}{\delta}}.
\end{align*}
Due to our assumption for a large value of $F_2(T_6)$, the function $F_2=F_2(t)$ blows up in finite time.

\subsection{Finite time blow-up of $F_1=F_1(t)$}\label{Subsection F1}
Let us multiply \eqref{1EQ 3} by $\frac{dF_1}{dt}(t)+F_1(t)$ to get
\begin{align*}
\frac{d^2F_2}{dt^2}(t)\,\Big(\frac{dF_1}{dt}(t)+F_1(t)\Big)\geq k_4e^{-\beta_3t}(t+R)^{-\beta_2}(F_1(t))^q\Big(\frac{dF_1}{dt}(t)+F_1(t)\Big).
\end{align*}
 By using integration by parts we derive
\begin{align*}
\int_{T_3}^t\frac{d^2F_2}{d\tau^2}(\tau)\,\Big(\frac{dF_1}{d\tau}(\tau)+F_1(\tau)\Big)\,d\tau&=\frac{dF_2}{dt}(t)\,\Big(\frac{dF_1}{dt}(t)+F_1(t)\Big)
-\frac{dF_2}{dt}(T_3)\,\Big(\frac{dF_1}{dt}(T_3)+F_1(T_3)\Big)\\
&\quad\,\,-\int_{T_3}^t\frac{dF_2}{d\tau}(\tau)\,\Big(\frac{d^2F_1}{d\tau^2}(\tau)+\frac{dF_1}{d\tau}(\tau)\Big)\,d\tau\\
&\leq \frac{dF_2}{dt}(t)\,\Big(\frac{dF_1}{dt}(t)+F_1(t)\Big).
\end{align*}
It is clear that for $t\geq T_7>T_3$ we have
\begin{align*}
&\int_{T_3}^t(\tau+R)^{-\beta_2}(F_1(\tau))^q\Big(\frac{dF_1}{d\tau}(\tau)+F_1(\tau)\Big)\,d\tau\\
&\geq\frac{1}{q+1}(t+R)^{-\beta_2}\big((F_1(t))^{q+1}-(F_1(T_3))^{q+1}\big)+\int_{T_3}^t(\tau+R)^{-\beta_2}(F_1(\tau))^{q+1}\,d\tau\\
&\geq\frac{1}{2(q+1)}(t+R)^{-\beta_2}(F_1(t))^{q+1},
\end{align*}
where we used the monotonically increasing property of $F_1=F_1(t)$. Hence,
\begin{align}\label{EQ 9}
\frac{dF_2}{dt}(t)\,\Big(\frac{dF_1}{dt}(t)+F_1(t)\Big)\geq \frac{k_4}{2(q+1)}e^{-\beta_3t}(t+R)^{-\beta_2}(F_1(t))^{q+1}.
\end{align}
Due to the fact that
\begin{align*}
(F_1(t))^{q+1}\Big(\frac{dF_1}{dt}(t)+F_1(t)\Big)\geq\frac{1}{q+2} \frac{d}{dt}\big((F_1(t))^{q+2}\big),
\end{align*}
we may multiply \eqref{EQ 9} by $\frac{dF_1}{dt}(t)+F_1(t)$ and integrate it over $[T_7,t]$ to obtain
\begin{align}\label{EQ 10}
F_2(t)\geq\frac{k_4}{4(q+1)(q+2)} \Big(\frac{dF_1}{dt}(t)+F_1(t)\Big)^{-2}(t+R)^{-\beta_2}(F_1(t))^{q+2},
\end{align}
where we used
\begin{align*}
\int_{T_7}^t\frac{dF_2}{d\tau}(\tau)\,\Big(\frac{dF_1}{d\tau}(\tau)+F_1(\tau)\Big)^2\,d\tau&
=F_2(t)\Big(\frac{dF_1}{dt}(t)+F_1(t)\Big)^2-F_2(T_7)\Big(\frac{dF_1}{dt}(T_7)+F_1(T_7)\Big)^2\\
&\quad\,\,-2\int_{T_7}^tF_2(\tau)\Big(\frac{dF_1}{d\tau}(\tau)+F_1(\tau)\Big)\Big(\frac{d^2F_1}{d\tau^2}(\tau)+\frac{dF_1}{d\tau}(\tau)\Big)\,d\tau\\
&\leq F_2(t)\Big(\frac{dF_1}{dt}(t)+F_1(t)\Big)^2,
\end{align*}
and
\begin{align*}
\int_{T_7}^t(\tau+R)^{-\beta_2}\frac{d}{d\tau}\big((F_1(\tau))^{q+2}\big)\,d\tau\geq\frac{1}{2} e^{-\beta_3t}(t+R)^{-\beta_2}(F_1(t))^{q+2},
\end{align*}
for $t\geq T_8>T_7$. Furthermore, plugging \eqref{EQ 10} into \eqref{1EQ 1} yields
\begin{align*}
\frac{1}{2p+1}\frac{d}{dt}\Big(\Big(\frac{dF_1}{dt}(t)+F_1(t)\Big)^{2p+1}\Big)&=\frac{d}{dt}\Big(\frac{dF_1}{dt}(t)+F_1(t)\Big)\,\Big(\frac{dF_1}{dt}(t)+F_1(t)\Big)^{2p}\\
&\geq\frac{k_2k_4^p}{(4(q+1)(q+2))^p}e^{-\beta_3t}(t+R)^{-\alpha_2-\beta_2p}(F_1(t))^{(q+2)p}.
\end{align*}
Finally, we multiply the above inequality by $\frac{dF_1}{dt}(t)+F_1(t)$ and integrate it over $[T_8,t]$ to obtain
\begin{equation}\label{F_1(t)}
\begin{aligned}
\frac{dF_1}{dt}(t)+F_1(t)\geq&\frac{k_2^{\frac{1}{2(p+1)}}k_4^{\frac{p}{2(p+1)}}}{2((q+2)p+1)^{\frac{1}{2(p+1)}}(4(q+1)(q+2))^{\frac{p}{2(p+1)}}}e^{-\frac{\beta_3}{2(p+1)}t}(t+R)^{-\frac{\alpha_2+\beta_2p}{2(p+1)}}(F_1(t))^{\frac{(q+2)p+1}{2(p+1)}}\\
&=:k_7e^{-\frac{\beta_3}{2(p+1)}t}(t+R)^{-\frac{\alpha_2+\beta_2p}{2(p+1)}}(F_1(t))^{\frac{(q+2)p+1}{2(p+1)}},
\end{aligned}
\end{equation}
where we used for $t\geq T_9>T_8$
\begin{align*}
&\int_{T_8}^t\Big(\frac{dF_1}{d\tau}(\tau)+F_1(\tau)\Big)\frac{d}{d\tau}\Big(\Big(\frac{dF_1}{d\tau}(\tau)+F_1(\tau)\Big)^{2p+1}\Big)\,d\tau\\
&\qquad =\Big(\frac{dF_1}{dt}(t)+F_1(t)\Big)^{2(p+1)}-\Big(\frac{dF_1}{dt}(T_8)+F_1(T_8)\Big)^{2(p+1)}\\
&\qquad \quad\,\,-\int_{T_8}^t\Big(\frac{d^2F_1}{d\tau^2}(\tau)+\frac{dF_1}{d\tau}(\tau)\Big)\Big(\frac{dF_1}{d\tau}(\tau)+F_1(\tau)\Big)^{2p+1}\,d\tau\\
&\qquad \leq\Big(\frac{dF_1}{dt}(t)+F_1(t)\Big)^{2(p+1)}.
\end{align*}
To prove our desired lemma, we shall distinguish between two cases.
\subsubsection{The condition $\alpha_2+\beta_2p<\alpha_1(pq-1)+2(p+1)$ holds.}
\noindent  Applying our derived estimate \eqref{1EQ 2} we have
\begin{align*}
\frac{dF_1}{dt}(t)+F_1(t)\geq k_0^{\epsilon}k_7e^{-\frac{\beta_3}{2(p+1)}t}(t+R)^{-\frac{\alpha_2+\beta_2p}{2(p+1)}+\alpha_1\epsilon}(F_1(t))^{\frac{(q+2)p+1}{2(p+1)}-\epsilon},
\end{align*}
with a positive constant $\epsilon>0$ to be determined later. From our assumption, we know
\begin{align*}
\left(\frac{\alpha_2+\beta_2p}{2(p+1)}-1\right)\frac{1}{\alpha_1}<\epsilon<\frac{(q+2)p+1}{2(p+1)}-1,
\end{align*}
which implies
\begin{align*}
-\frac{\alpha_2+\beta_2p}{2(p+1)}+\alpha_1\epsilon>-1\quad\text{and}\quad
\frac{(q+2)p+1}{2(p+1)}-\epsilon>1.
\end{align*}
Let us consider the Cauchy problem
\begin{align} \label{auxiliaryCauchyproblem}
\frac{dZ}{dt}(t)+Z(t)=\kappa e^{-\gamma t}(t+R)^{-\alpha}(Z(t))^{\beta},\quad Z(T_9)=F_1(T_9),
\end{align}
that is an auxiliary initial value problem to \eqref{F_1(t)}.
The solution to \eqref{auxiliaryCauchyproblem} is explicitly given by
\begin{align*}
Z(t)&=e^{- t}\Big((Z(T_9))^{1-\beta}-(\beta-1)\nu\int_{T_9}^te^{-(\beta+\gamma-1)\tau}(\tau+R)^{-\alpha}\,d\tau\Big)^{-\frac{1}{\beta-1}}.
\end{align*}
Thus, if $-\alpha>-1$ and $\beta>1$, then the solution $Z=Z(t)$ blows up when $Z(T_9)$ is large. According to the Petrovitsch theorem, we conclude that $F_1=F_1(t)$ blows up in finite time, too.
\begin{rem} If we drop the assumption on large value of $F_1(t_0)$ in the theorem, this part cannot be completed due to \eqref{F_1(t)}. It seems to be difficult to avoid $e^{-t}$ in the last step of the proof.\end{rem}
\subsubsection{The condition $\alpha_2+\beta_2p=\alpha_1(pq-1)+2(p+1)$ holds.}
\noindent Similar as in Section \ref{Sec1000}, we multiply $(F_1(t))^{-1-\delta}$ with $\delta\in\big(0,\big(\frac{\alpha_2+\beta_2p}{2(p+1)}-1\big)\frac{1}{\alpha_1}\big)$ to \eqref{F_1(t)} and apply the estimate \eqref{1EQ 2} to get
\begin{align*}
-\frac{1}{\delta}\frac{d}{dt}\big((F_1(t))^{-\delta}\big)+(F_1(t))^{-\delta}&\geq k_7e^{-\frac{\beta_3}{2(p+1)}t}(t+R)^{-\frac{\alpha_2+\beta_2p}{2(p+1)}}(F_1(t))^{\frac{(q+2)p+1}{2(p+1)}-1-\delta}\\
&\geq k_7e^{-\frac{\beta_3}{2(p+1)}t}(t+R)^{-\frac{\alpha_2+\beta_2p}{2(p+1)}}(F_1(t))^{\big(\frac{\alpha_2+\beta_2p}{2(p+1)}-1\big)\frac{1}{\alpha_1}-\delta}\\
&\geq k_7e^{-\frac{\beta_3}{2(p+1)}t}(t+R)^{-1-\delta\alpha_1},
\end{align*}
where we used our condition
\begin{align*}
\Big(\frac{\alpha_2+\beta_2p}{2(p+1)}-1\Big)\frac{1}{\alpha_1}=\frac{(q+2)p+1}{2(p+1)}-1.
\end{align*}
Thus, by direct computation we have
\begin{align*}
F_1(t)&\geq e^{T_9-t}\Big((F_1(T_9))^{-\delta}-\delta k_7e^{\delta T_9}\int_{T_9}^te^{-\delta\tau-\frac{\beta_3}{2(p+1)}\tau}(\tau+R)^{-1-\delta\alpha_1}d\tau\Big)^{-\frac{1}{\delta}}\\
&\geq e^{T_9-t}\Big((F_1(T_9))^{-\delta}-\frac{k_7}{\alpha_1}e^{\delta T_9-\delta t}\big((T_9+R)^{-\delta\alpha_1}-(t+R)^{-\delta\alpha_1}\big)\Big)^{-\frac{1}{\delta}}.
\end{align*}
Due to our assumption for a large value of $F_1(T_9)$, the function $F_1=F_1(t)$ blows up in finite time.
\section{Proof of our main results}\label{Section proof main result}
\subsection{Proof of Theorem \ref{Local existence of solution}}\label{Subsection proof of local solution}
For the cases $n=1,2$, one can see \cite{Strauss1989,Wakasugi2017}. We may apply the Gagliardo-Nirenberg inequality and Banach's fixed-point theorem to prove that there exists a uniquely determined solution $(u^*,v^*)\in X(T)\times X(T)$ for a positive $T$, where the evolution space $X(T)$ is defined by
\begin{equation*}
X(T):=\big\{f\in\ml{C}\big([0,T],H^1(\mb{R}^n)\big)\cap \ml{C}^1\big([0,T],L^2(\mb{R}^n)\big):\,\,\text{supp }f(t,\cdot)\subset\{|x|\leq t+R\}\big\}.
\end{equation*}
Here we only need to restrict to $1<p,q<\infty$.\\
For the remaining cases $n\geq3$, one can combine the proofs stated in \cite{ReissigEbert2018,Sideris1984,Strauss1989,Wakasugi2017}. We may apply the embedding theorem, the Gagliardo-Nirenberg inequality and Banach's fixed-point theorem to prove that there exists a uniquely determined solution $(u^*,v^*)\in X(T)\times Y(T)$ for a positive $T$, where the evolution space $Y(T)$ is defined by
\begin{equation*}
Y(T):=\big\{f\in\ml{C}\big([0,T],L^{\frac{2(n+1)}{n-1}}(\mb{R}^n)\big):\,\,\text{supp }f(t,\cdot)\subset\{|x|\leq t+R\}\big\}.
\end{equation*}
On the one hand, the restriction of the exponent $p$ to $1<p\leq\frac{n+3}{n-1}$ comes from the estimate
\begin{align*}
\big\||v(t,\cdot)|^p\big\|_{L^{\frac{2(n+1)}{n+3}}(\mb{R}^n)}\leq c\|v\|^p_{Y(T)}\quad\text{if}\quad\frac{2(n+1)p}{n+3}\leq \frac{2(n+1)}{n-1},
\end{align*}
with $c>0$. On the other hand, the restriction of the exponent $q$ to $1<q\leq\frac{n}{n-2}$ for $n\geq3$ comes from the application of the Gagliardo-Nirenberg inequality.
\subsection{Proof of Theorem \ref{Blow-up result}}\label{Subsection Poof blowup}
Let us choose
\begin{align*}
\alpha_1&=1+\frac{2-p}{2}(n-1),\qquad\alpha_2=n(p-1),\\
\beta_1&=1,\qquad \beta_2=n(q-1),
\end{align*}
in Lemma \ref{Generalzied Kato lemma}, we immediately get from \eqref{Cond 1} and \eqref{Cond 2}, respectively,
\begin{align*}
\frac{q+1}{pq-1}\leq\frac{n-1}{2}\quad\text{and}\quad\frac{2+2p^{-1}}{pq-1}\leq\frac{n-1}{2}.
\end{align*}
Then, the proof of Theorem \ref{Blow-up result} is completed.

\section*{Acknowledgments}
The PhD study of Mr. Wenhui Chen are supported by S\"achsisches Landesgraduiertenstipendium.

\bibliographystyle{elsarticle-num}

\end{document}